\documentclass[12pt]{article}
\usepackage{amsmath}
\usepackage{amssymb}

\newcommand{\co}{\mathbb{C}}
\newcommand{\na}{\mathbb{N}}

\newcommand{\rl}{\mathbb{R}}

\newcommand{\vep}{\varepsilon}
\newcommand{\dd}{\mathbb{D}}
\newcommand{\tz}{\mathbb{T}}
\newcommand{\cp}{\mathbb{C}^+}

\newcommand{\kp}{{K^p_{\Theta}}}
\newcommand{\kd}{{K^2_{\Theta}}}

\title{Compact embeddings of model subspaces \\ of the Hardy space}

\author {Anton D. Baranov}
\date{}

\begin {document}
\maketitle
\sloppy
\noindent
{\small {\bf Abstract.}
We study embeddings of model (star-invariant) 
subspaces $K^p_{\Theta}$ of the Hardy space $H^p$, associated
with an inner function $\Theta$. We obtain a criterion  for 
the compactness of the embedding of $K^p_{\Theta}$ into $L^p(\mu)$
analogous to the Volberg--Treil theorem on bounded embeddings
and answer a question posed by Cima and Matheson. 
The proof is based on Bernstein inequalities for functions in $K^p_{\Theta}$.
Also we study measures $\mu$ such that the 
embedding operator belongs to a Schatten--von Neumann ideal.
}
\\
\medskip

\begin{center}
{\bf 1. Introduction and main results}
\end{center}

Let $\dd =\{z:\,|z|<1\}$ be the unit disc and let $\tz=\{z:\,|z|=1\}$
be the unit circle. We denote by $m$ the normalized Lebesgue measure on $\tz$.
A function $\Theta$ which is  analytic and bounded in $\dd$
is said to be {\it inner} if
$|\Theta| =1$ $m$-a.e. on $\tz$ in the sense of nontangential boundary values.
Recall that each inner
function $\Theta$ admits the factorization (up to a constant
unimodular factor)
$$
\Theta(z)= B(z)I_\psi(z)
$$
where 
$$
B(z)=\prod\limits_n \frac{|z_n|}{z_n}\cdot
\frac{z_n-z}{1- \overline z_n z}, \qquad z\in\dd,
$$
is the {\it Blaschke product} 
with zeros $z_n\in\mathbb{D}$, and the sequence $\{z_n\}$ 
satisfies the Blaschke condition $\sum_n (1-|z_n|)<\infty$
(we assume $|z_n|/z_n = 1$ if $z_n=0$). 
The {\it singular inner function} $I_\psi$ is defined by the formula
$$
I_\psi(z)=\exp\bigg(- \int\limits_\mathbb{T}\frac{\zeta+z}
{\zeta-z}\, d\psi(\zeta)\bigg), \qquad z\in\dd,
$$
where $\psi$ is a finite Borel measure on $\tz$
singular with respect to $m$.
\smallskip

Let $H^p$ denote the Hardy space in $\dd$, $1\le p\le \infty$. 
With each inner function $\Theta$ we associate the subspace
$$
K^p_{\Theta}=H^p\cap \Theta \overline{H^p_0}
$$
where $H_0^p = \{f\in H^p: f(0)=0\}= zH^p$. Equivalently, 
one can define $K_\Theta^p$ as the set of all functions $f$ in $H^p$
such that $\langle f,\Theta g \rangle=
\int_\tz f\, \overline{\Theta g} \,  dm=0$
for any $g\in H^q$, $1/p+1/q =1$. Note that $K_\Theta^2 = H^2\ominus 
\Theta H^2$. It is well known that, for $1\le p<\infty$, any closed subspace 
of $H^p$ invariant with respect to the backward shift 
$(S^*f)(z)=\frac{f(z)-f(0)}{z}$
is of the form $K_\Theta^p$ for some inner function $\Theta$ 
(see \cite[Chapter II]{ga} and \cite{al9}). Subspaces 
$K_\Theta^p$ are often called {\it star-invariant subspaces}. These subspaces 
play an outstanding role both in function and operator theory
(see \cite{cimros, hnp, nik, nk2}) and, in particular,
in the Nagy--Foias model for contractions in a Hilbert space
(therefore they are sometimes refered to  as {\it model subspaces}).
Note that if $\Theta$ is a Blaschke product, then $\kp$
coincides with the closed linear span of simple fractions
with the poles of corresponding multiplicities at the points 
$1/\overline z_n$.
\smallskip

Let $\sigma(\Theta)$ be the so-called {\it spectrum} of the inner
function $\Theta$, that is, the set of all 
$\zeta\in\overline {\mathbb{D}}$
such that $\liminf\limits_{z\to\zeta,\, z\in\mathbb{D}}|\Theta(z)|=0$.
Equivalently, $\sigma(\Theta)$ is the smallest closed subset of
$\overline {\mathbb{D}}$ containing zeros $z_n$
and the support of the measure $\psi$.
Clearly, $\Theta$, as well as any element of $\kp$, 
has an analytic extension across any subarc  
of the set $\mathbb{T}\setminus \sigma(\Theta)$.
\smallskip

In the present paper we study the following problem:
given an inner function $\Theta$ and $p\ge 1$, 
describe the class of Borel measures $\mu$ in the 
closed disc $\overline{\dd}$ such that the space
${K^p_{\Theta}}$ is embedded into $L^p(\mu)$
or such that the embedding is compact.
This problem was posed by Cohn in 1982 \cite{con1}: 
in spite of a number of partial results it is still open. 
The embedding ${K^p_{\Theta}}\subset L^p(\mu)$
is equivalent to the estimate
\begin{equation}
\label{emb}
\|f\|_{L^p(\mu)}\le C\|f\|_p,\qquad
f\in{K^p_{\Theta}}.
\end{equation}
The class of such measures $\mu$ is denoted by ${\cal C}_p(\Theta)$.
\smallskip

Recall that a finite Borel measure $\mu$ in the closed unit 
disc $\overline{\dd}$
is said to be a {\it Carleson measure} if there is
a constant $M>0$ such that 
\begin{equation}
\label{carl}
\mu(S(I))\le M|I|
\end{equation}
for any arc  $I\subset \mathbb{T}$. 
Here and in  what follows we denote by $|I|$ the length of 
an arc $I$, and by $S(I)$ the {\it Carleson square}
\begin{equation}
\label{csq}
S(I)=\{z = \rho e^{i\varphi} \in \overline{\dd}:\, 
e^{i\varphi}\in I,\;  1-(2\pi)^{-1}|I| \le \rho \le 1\}.
\end{equation}
We denote the class of Carleson measures by ${\cal C}$;
for a measure $\mu\in {\cal C}$ we denote by $M_\mu$ the smallest constant
$M$ in (\ref{carl}). The classical Carleson theorem states that
$H^p\subset L^p(\mu)$ for some (any) $p>0$
if and only if $\mu \in {\cal C}$.
The embedding $H^p\subset L^p(\mu)$ is compact if and only if 
$\mu$ is a {\it vanishing Carleson measure}, that is,
\begin{equation}
\label{csq1}
\lim\limits_{|I|\to 0}\frac{\mu(S(I))}{|I|} = 0
\end{equation}
(see \cite{power}; some generalizations can be found in \cite{blasco}).

Obviously, ${\cal C} \subset {\cal C}_p(\Theta)$. One may expect that
the class ${\cal C}_p(\Theta)$ will depend essentially on geometric
properties of $\Theta$. At present the class 
${\cal C}_p(\Theta)$ is described explicitly only for some very
special classes of inner functions. We say
that an inner function satisfies the {\it connected level set condition}
(we write $\Theta\in CLS$) if the level set 
$$
\Omega(\Theta,\varepsilon) = \{z\in \dd: |\Theta(z)|<\vep\}
$$
is connected for some $\vep\in (0,1)$.
Cohn \cite{con1} showed that if $\Theta \in CLS$, then
it suffices to verify inequality (\ref{emb}) for
the reproducing kernels of the space ${K^2_{\Theta}}$ 
(see Section 2 for the definition).
Recently Nazarov and Volberg \cite{nv}
showed that this is no longer true in the general case.
\smallskip

A geometric condition on $\mu$ sufficient for the embedding
of ${K^p_{\Theta}}$ is due to  Volberg and Treil \cite{vt}:
the embedding ${K^p_{\Theta}}\subset L^p(\mu)$ takes place
if there is an $\varepsilon\in(0,1)$ such that
$\mu(S(I))\le C|I|$ for all squares $S(I)$
such that
$$
S(I)\cap\Omega(\Theta,\varepsilon)\ne\emptyset.
$$
Thus, we should check the Carleson condition
(\ref{carl}) only for squares of the special form.
Denote by ${\cal C}(\Theta)$ the class of measures satisfying the conditions of
the Volberg--Treil theorem for some $\varepsilon\in (0,1)$.
It was proved by Aleksandrov \cite{al2}
that the condition $\mu\in{\cal C}(\Theta)$ is necessary
(that is ${\cal C}_p(\Theta)={\cal C}(\Theta)$) if and only if 
$\Theta \in CLS$. Moreover, if $\Theta$ is not a $CLS$ function, then the class 
${\cal C}_p(\Theta)$, $p>0$, depends essentially on
the exponent $p$ (in contrast to the classical Carleson theorem). 
Some other embedding theorems
were obtained in \cite{al2,con2,d5}, compact
embeddings are considered in \cite{cimmat,con2,v81}.
\smallskip

Of special interest is the case
when $\mu=\sum_{n\in \na} a_n \delta_{\lambda_n}$ is a discrete measure;
then embedding is equivalent to the Bessel property for the system of 
reproducing kernels $\{k_{\lambda_n}\}$.
Also the particular case when $\mu$ is a measure on
the unit circle is of great importance. In contrast to the embeddings of
the whole Hardy class $H^p$ (note that Carleson measures on $\tz$
are measures with bounded density with respect to Lebesgue measure $m$),
the class ${\cal C}_p(\Theta)$ always contains nontrivial examples  of singular
measures on $\tz$; in particular, for $p=2$, the Clark measures \cite{cl}
for which the embedding $K^2_{\Theta}\subset L^2(\mu)$ is isometric. 
On the other hand, if $\mu=wm$, $w\in L^2(\tz)$,
then the embedding problem is related to the 
properties of the Toeplitz operator $T_w$ \smallskip
(see \cite{con2}).

A new approach to the embedding theorems was suggested by the author
in \cite{bar03, bar05}. It is based on Bernstein inequalities
for $\kp$. By a Bernstein inequality we mean an estimate 
for a weighted norm of the derivative $f'$ 
in terms of the standard $L^p$-norm of $f\in K_\Theta^p$, that is,
an estimate of the form
\begin{equation}
\label{wei}
\|f'\|_{L^p(\mu)}\le C\|f\|_p, \qquad f\in{K^p_{\Theta}},
\end{equation}
where $\mu$ is a measure in the closed disc $\overline{\dd}$.
This approach made possible to obtain essentially new 
embedding theorems generalizing the Volberg--Treil theorem 
as well as certain theorems due to Cohn. One more application 
of Bernstein inequalities is connected with the stability
of Bessel sequences and Riesz bases of reproducing kernels
\cite{bar05a}. The results in \cite{bar05, bar05a}
are obtained in the context of the Hardy spaces in the half-plane.
In Section 3 we discuss their analogs for the disc.
\smallskip

In the present paper we apply Bernstein inequalities 
for $\kp$ to study compactness and Schatten--von Neumann class properties
of the embeddings. One of our main results 
is a geometric condition sufficient 
for the compactness of the embedding
which is similar to the Volberg--Treil Theorem.
Now we need to verify the "vanishing condition" (\ref{csq1}) 
only for squares intersecting the level sets.
This condition turns out to be also necessary in  the case 
of $CLS$ inner functions.
\medskip
\\
{\bf Theorem 1.1.}
{\it Let $1 < p<\infty$, let $\mu$ be a Borel measure on $\overline{\dd}$,
and let $\varepsilon \in (0,1)$. Then $(i)$ implies $(ii)$ where
\smallskip

$($i$)$ for each $\eta>0$ there exists a $\delta>0$ such that
$\mu(S(I))/|I| <\eta$
whenever $|I|<\delta$ and $S(I)\cap \Omega(\Theta,\varepsilon)\ne\emptyset$;
\smallskip

$($ii$)$ the subspace $K^p_{\Theta}$ is compactly embedded in $L^p(\mu)$.
\smallskip
\\
Conversly, if $\Theta \in CLS$, then $($ii$)$ implies $($i$)$.}
\medskip

Implication ($ii$)$\Longrightarrow$($i$) for $CLS$ inner functions
was proved by Cima and Matheson in \cite{cimmat}, where 
the question was posed whether the converse is true. 
Theorem 1.1 provides a positive answer.
Also in \cite{cimmat} another vanishing condition for a measure $\mu$
was introduced which is sufficient for the compactness of the embedding
$K^p_{\Theta} \subset L^p(\mu)$ for any $p>0$.
We show (see Proposition 4.1) that this condition implies 
condition ($i$) of Theorem 1.1 (thus, we answer 
one more question posed in \cite{cimmat}).
Theorem 1.1 is deduced from a more general embedding theorem (Theorem 3.1)
which is an extension of the Volberg--Treil  theorem.
\smallskip

In Sections 5-6 conditions ensuring the inclusion of the embedding 
operator ${\cal J}_\mu: {K^2_\Theta}\to L^2(\mu)$, ${\cal J}_\mu f=f$, 
into the Schatten--von Neumann ideal 
${\cal S}_r$ are studied. We give a complete description of such measures 
for the case $\Theta \in CLS$ and $r\ge 1$. 
For $\varepsilon \in (0,1)$ consider 
a Whitney-type decomposition of the set 
$\tz \setminus \sigma(\Theta)$ into the union of arcs $I_k$
with the property 
$$
dist\, (I_k, \Omega(\Theta,\varepsilon)) \asymp |I_k|
$$
(see Section 3, Lemma 3.5, for the details).
\medskip
\\
{\bf Theorem 1.2.}
{\it Let $\mu$ be a Borel measure with $supp\,\mu \subset 
\bigcup\limits_k S(I_k)$. Assume that for some $r>0$,
\begin{equation}
\label{23-a}
\sum\limits_{k} \left(\frac{\mu(S(I_k)}{|I_k|}\right)^{r/2}<\infty.
\end{equation}
Then ${\cal J}_\mu \in {\cal S}_r$. }
\medskip

We denote by $R_{n,m}$ the elements of the standard 
dyadic partition of the disc (see the definition in Section 5).
Then we have the following necessary condition.
\medskip
\\
{\bf Theorem 1.3.}
{\it Let $\varepsilon\in (0,1)$. If ${\cal J}_\mu
\in {\cal S}_r$, $r\ge 1$, then
\begin{equation}   
\label{dya2}
\sum\limits_{R_{n,m}\cap \Omega(\Theta, \varepsilon) \ne\emptyset} 
(2^{n}\mu(R_{n,m}))^{r/2}<\infty.
\end{equation}  }

For $CLS$ inner functions the statements converse to 
Theorems 1.2 and 1.3 are also true; 
in this case we give a complete description of embeddings 
of a class ${\cal S}_r$, 
\medskip
$r\ge 1$. 
\\
{\bf Theorem 1.4.}
{\it Let $\Theta \in CLS$, let $\vep \in (0,1)$, and let 
$\mu$ be a Borel measure in $\overline{\mathbb{D}}$. 
Then the embedding operator ${\cal J}_\mu$ 
belongs to ${\cal S}_r$, $r\ge 1$, if and only if
$\mu$ satisfies (\ref{23-a}) and 
\medskip
(\ref{dya2}). }

Our conditions, stated in terms of dyadic partition 
of the disc, are similar to a theorem of Luecking 
\cite{luek} which characterizes Schatten--von Neumann class 
embeddings of the whole Hardy class, as well as to 
results of Parfenov \cite{parf1, parf2}.
For further discussion see Sections 5-6.
\smallskip

As in \cite{bar05, bar05a}, our main tools 
in the present paper are Bernstein inequalities for 
$\kp$. For the sake of completeness we include a discussion of these
results which are of independent interest. Moreover, in contrast to
\cite{bar05}, we give estimates for higher order derivatives as well.
Now we state two results of this type. The first of them shows that 
the growth of the derivative is controlled by the distance
to a level set. For $\zeta\in\mathbb{T}$ let
$d_\varepsilon(\zeta)=dist\, (\zeta,\Omega(\Theta,\varepsilon))$.
\medskip
\\
{\bf Theorem 1.5.}
{\it Let $\varepsilon\in (0,1)$,  $1<p<\infty$. Then 
$$
\|f^{(n)}\cdot d^n_\varepsilon\|_p \le C(p,n,\varepsilon)
\|f\|_p, \qquad f\in {K^p_{\Theta}}.
$$}

Let $\Theta \in CLS$. Then the boundary spectrum 
$\sigma(\Theta)\cap\tz$ has zero Lebesgue measure \cite{al1};
thus, the $n$th derivative $f^{(n)}(\zeta)$ is well-defined
for almost all $\zeta\in\mathbb{T}$. 
\medskip
\\
{\bf Theorem 1.6.}
{\it Let $\Theta \in CLS$, $1<p<\infty$, $\mu \in {\cal C}$. Then
\begin{equation}
\label{ogr0}
\int\limits_\dd \left|\frac{f^{(n)}(z)}{\|k_z\|_2^{2n}}\right|^p 
d\mu(z) \le C(\Theta, p, n, \mu) \|f\|_p^p,
\qquad f\in {K^p_{\Theta}},
\end{equation}
where $\|k_z\|_2$ denotes the $L^2$-norm of the reproducing kernel of $\kd$.
In particular, 
\begin{equation}
\label{lev-cls}
\|f^{(n)} \cdot |\Theta'|^{-n} \|_p\le C(\Theta,p, n)
\|f\|_p, \qquad f\in{K^p_{\Theta}}.
\end{equation}}
\smallskip

Theorems 1.5 and 1.6 are corollaries of a much more general,
but somewhat more complicated Bernstein inequality which
will be proved in Section 2 (Theorem 2.1).
\medskip

We make use of the following notations:
given nonnegative functions $g$ and $h$, we write
$g \lesssim h$ if $g\le Ch$ for a positive constant $C$
and all admissible values of the variables;
we write  $g\asymp h$ if $g\lesssim h\lesssim g$.
Letters $C$, $C_1$, etc. will denote various positive constants
which may change their values in different occurrences.
\bigskip

\begin{center}
{\bf 2. Bernstein inequalities for higher order derivatives}
\end{center}

From now on we assume that $p\in [1,\infty)$ and $q$ 
is the conjugate exponent, that is, $1/p+1/q=1$;
by $L^p$ we denote the standard space $L^p(\mathbb{T}, m)$.
\smallskip

We start with a discussion of the local behavior of the
elements of model subspaces and their derivatives near the boundary.
In what follows, the reproducing kernels 
play the most prominent role. The function 
$$
k_z(\zeta)=
\frac{1-\overline{\Theta(z)}\Theta(\zeta)}{1-\overline z \zeta}
$$
is the reproducing kernel of the space ${K^2_{\Theta}}$
corresponding to a point $z\in\mathbb{D}$.
Since $k_z \in K^\infty_\Theta$, we have
\begin{equation}
\label{inte0}
f(z)=\int\limits_\mathbb{T} f(\tau)\overline{k_z(\tau)} \, dm(\tau), 
\qquad f\in{K^p_{\Theta}},
\end{equation}
for any  $1\le p \le \infty$. We have an analogous representation
for the $n$th derivative,
\begin{equation}
\label{inte2}
f^{(n)}(z)=n! \int\limits_\mathbb{T} \overline\tau ^{n}
f(\tau) (\overline{k_z(\tau)})^{n+1} \, dm(\tau),\qquad f\in{K^p_{\Theta}},
\end{equation}
which follows from the fact that $(1-\tau\overline z)^{n+1}-
(k_z(\tau))^{n+1} \in \Theta H^\infty$.
\smallskip

Integral representations (\ref{inte0})--(\ref{inte2}) 
may be extended to certain points $z=\zeta$ on the unit circle $\tz$.
It is well known that any element of $K^p_{\Theta}$ has an 
analytic extension across any subarc  
of the set $\mathbb{T}\setminus \sigma(\Theta)$, and thus 
(\ref{inte0})--(\ref{inte2})
hold for any $\zeta \in \mathbb{T}\setminus \sigma(\Theta)$. 
The boundary behavior 
at the points of the spectrum is more subtle and depends 
on the ``density'' of the spectrum
near the given point. For $\zeta\in \tz$, put
$$
S_q(\zeta)=\sum\limits_n \frac{1-|z_n|^2}{|\zeta- z_n |^{q}}+
\int\limits_\mathbb{T}\frac{d\psi(\tau)}{|\zeta-\tau|^{q}}.
$$
Then, by the results of Ahern and Clark \cite{ac70} 
and Cohn \cite{cohn86}, $f^{(n)}(\zeta)$ (understood in the sense 
of nontangential boundary values) exists for any 
$f\in K^p_\Theta$ if and only if $S_{(n+1)q}(\zeta)<\infty$;
in this case $k^{n+1}_\zeta \in K^q_\Theta$ and (\ref{inte2}) 
holds with $z=\zeta$.
The quantity $S_2$, which is responsible for the inclusion $k_\zeta \in \kd$,
is of special importance. In the case when $\Theta$ has a nontagential limit  
at $\zeta$ and $\Theta(\zeta) \in \tz$, it coincides 
with the modulus of the nontangential derivative
of $\Theta$ at $\zeta$ (by the nontangential derivative we mean
$\lim\limits_{z\to \zeta}\frac{\Theta(z)-\Theta(\zeta)}
{z-\zeta}$ where $z$ tends to $\zeta$ nontangentially):
$$
|\Theta'(\zeta)|=\sum\limits_n \frac{1-|z_n|^2}{|\zeta-z_n |^2}+
\int\limits_\mathbb{T}\frac{d\psi(\tau)}{|\zeta - \tau|^2}.
$$

Our main result in this section is the following 
weighted Bernstein inequality of the form (\ref{wei})
which holds for an arbitrary inner function $\Theta$ and for the measures of 
the form $w\mu$ where $\mu$ is a Carleson
measure and the weight $w(z)$ depends on the norm
of the kernel $k^{n+1}_z$ in $L^q$ (that is, in essence,
on the norm of the functional $f\mapsto f^{(n)}(z)$, $f\in{K^p_{\Theta}}$).
Put
$$
w_{p,n} (z)=\|k^{n+1}_z\|_q^{-\frac{pn}{pn+1}}.
$$
We assume $\|k^{n+1}_\zeta\|_q=\infty$ and $w_{p,n}(\zeta)=0$
whenever $\zeta \in\mathbb{T}$ and $S_{(n+1)q}(\zeta)=\infty$;
thus $f^{(n)}(z)w_{p,n}(z)$ is well-defined for any $f\in {K^p_{\Theta}}$,
$z\in \overline{\dd}$.
\bigskip
\\
{\bf Theorem 2.1.} {\it Let $\mu\in {\cal C}$, $1\le p<\infty$.
Then the operator  
$$
(T_{p,n} f)(z)=f^{(n)}(z) w_{p,n}(z)
$$
is of weak type $(p,p)$ as an operator from ${K^p_{\Theta}}$ to $L^p(\mu)$
and is bounded as an operator from $K^r_{\Theta}$ to $L^r(\mu)$ 
for any $r>p$; moreover there is a constant $C=C(M_\mu,p,r,n)$ such that
\begin{equation}
\label{main}
\|f^{(n)}w_{p,n}\|_{L^r(\mu)}
\le C\|f\|_r,\qquad f\in{K^r_{\Theta}}.
\end{equation} }

To apply Theorem 2.1, one should have
effective estimates of the considered weights, that is,
of the norms of reproducing kernels. For $p=2$ we have an explicit formula:
$\|k_z\|_2^2=\frac{1-|\Theta(z)|^2}{1-|z|^2}$, $z\in\dd$, 
and $\|k_\zeta\|_2^2=|\Theta'(\zeta)|$, $\zeta\in\mathbb{T}$.
For $\Theta\in CLS$, sharp estimates for the norms are known 
(see (\ref{onecomp}) below).
We also relate the weight $w_{p,n}$ to
geometric properties of the level sets of the function $\Theta$:
\begin{equation}
\label{geomineq}
d_\vep(\zeta) \lesssim w_{p,n}(\zeta)
\lesssim |\Theta'(\zeta)|^{-1}, \qquad \zeta\in \mathbb{T}. 
\end{equation}
For the proof we refer to \cite[Lemma 4.5]{bar05} where the 
inequality is obtained for $n=1$; the arguments extend to general
$n$ in an obvious way. Some close results may be found in \cite{al2}.
\smallskip

The proof of Theorem 2.1 is based
on the integral representation (\ref{inte2}) which
reduces the study of differentiation operator
to the study of certain singular integral operators.
We deduce inequality (\ref{main}) from the boundedness
of the following integral operators in $L^p$-spaces 
associated with Carleson measures \cite[Theorems 3.1, 3.2]{bar05}. 
\medskip
\\
{\bf Theorem 2.2.} 
{\it Let $\mu\in{\cal C}$, and let $h$ be a nonnegative 
function in $\overline {\mathbb{D}}$,
measurable with respect to $\mu$ and $m$, and
such that $h(z)\ge A(1-|z|)$, $z\in\mathbb{D}$, 
for some constant $A>0$. Put
$$
Tf(z)=h(z)\int\limits_{|\zeta-z|\ge h(z)}\frac{f(\zeta)}{|\zeta-z|^2}\, 
dm(\zeta),\qquad z\in\overline{\mathbb{D}}.
$$
Then $T$ is of weak type $(1,1)$ as an operator from
$L^1$ to $L^1(\mu)$ and $T$ is a bounded operator 
from $L^p$ to $L^p(\mu)$ for any $p>1$.
Moreover, the norm of $T$ does not exceed some constant $C$ which depends
only on $p$, $A$ and the Carleson constant $M_\mu$ of the measure $\mu$.}
\medskip

We also consider a class of integral operators
with a ``diagonal'' kernel.
Let $\mu$, $\nu\in{\cal C}$ and let $K(z,u)$ be a 
$\mu\times\nu$-measurable function.
For $z\in\overline{\mathbb{D}}$ put
$$
\Delta_z(p)=\{u\in\overline{\mathbb{D}}: 
|u-z|<\|K(z,\cdot)\|_{L^q(\nu)}^{-p}\}
$$
(we assume without loss of generality that $K(z,\cdot)$ is $\nu$-measurable
for any $z$, and we put $\|K(z,\cdot)\|_{L^q(\nu)}^{-p}=0$ whenever
$K(z,\cdot)\notin L^q(\nu)$).
Consider the following ``truncation'' of the integral operator 
with the kernel $K$:
$$
T_pf(z)=\int\limits_{\Delta_z(p)}K(z,u)f(u)\,d\nu(u).
$$
{\bf Theorem 2.3.}
{\it If $\|K(z,\cdot)\|_{L^q(\nu)}^{-p}\ge A(1-|z|)$, then the operator 
$T_p$ is of weak type $(p,p)$ as an operator from
$L^p(\nu)$ to $L^p(\mu)$ and is bounded as an operator from 
$L^r(\nu)$ to  $L^r(\mu)$ for any $r>p$.}
\medskip

Detailed proofs of Theorems 2.2 and 2.3 are given in 
\cite{bar05} in the half-plane setting; the proofs
for the disc follow by exactly the same arguments
and we omit them.
\smallskip

The idea of the proof of Theorem 2.1 
is to split the integral, which represents the derivative,
into ``diagonal'' and ``off-diagonal''
parts and to estimate them separately making use of
Theorems 2.2 and 2.3, respectively.
\smallskip
\\
{\bf Proof of Theorem 2.1.}
Put $h(z) = \big(w_{p,n}(z)\big)^{1/n}$. A trivial inequality
$|k_z(w)|\le 2(1-|z|)^{-1}$, $z,w\in \dd$, implies $h(z)\ge A(1-|z|)$.
Multiply the integral in (\ref{inte2}) by $w_{p,n} (z)$ and split it
into two parts:
\begin{equation}
\label{c79}
\frac{1}{n!} w_{p,n} (z)f^{(n)}(z)= I_1f(z)+I_2f(z)
\end{equation}
where
$$
I_1f(z)= w_{p,n}(z)\int\limits_{|\zeta-z|\ge h(z)} 
\overline\zeta^n f(\zeta) 
\overline{k^{n+1}_z(\zeta)}\,dm(\zeta),
$$
$$
I_2f(z)= w_{p,n}(z)\int\limits_{|\zeta-z|< h(z)} \overline\zeta^n f(\zeta) 
\overline{k^{n+1}_z(\zeta)}\,dm(\zeta).
$$
Since $|1-\overline \zeta z|=|\zeta- z|$, $\zeta\in \tz$, we have
$$
|I_1f(z)|\le  C h^n(z)
\int\limits_{|\zeta-z|\ge h(z)}
\frac{|f(\zeta)|}{|\zeta -z|^{n+1}} \le C
h(z) \int\limits_{|\zeta-z|\ge h(z)}
\frac{|f(\zeta)|}{|\zeta -z|^2} \, dm(\zeta).
$$
Now Theorem 2.2 implies that the 
operator $I_1$ is bounded as an operator
from $L^r$ to $L^r(\mu)$ for any Carleson measure $\mu$
and $r>1$.

To estimate the integral $I_2f$, 
put $K(z,\zeta)=(h(z))^n \overline{k_z^{n+1}(\zeta)}$.
Then $\|K(z,\cdot)\|_q^{-p}=(h(z))^{-pn} \|k_z^{n+1} \|^{-p}_q=
h(z)$. Thus,
$$
I_2f(z)=\int\limits_{|\zeta-z|<\|K(z,\cdot)\|_q^{-p}} 
\overline\zeta^n f(\zeta) K(z,\zeta) \,dm(\zeta),
$$
and, applying Theorem 2.3, we conclude that
the operator $I_2$ is of weak type $(p,p)$ 
as an operator from $L^p$ to $L^p(\mu)$
and is bounded as an operator from $L^r$ to $L^r(\mu)$ for any 
$\mu\in{\cal C}$ and $r>p$. $\bigcirc$
\medskip
\\
{\bf Proof of Theorem 1.5.} In view of the inequality (\ref{geomineq}),
the statement follows from Theorem 2.1 with $\mu=m$. $\bigcirc$
\medskip
\\
{\bf Proof of Theorem 1.6.}
It is shown by Alexandrov \cite{al2} that, for $\Theta \in CLS$, we have
\begin{equation}
\label{onecomp}
\|k_z\|_s^s\asymp \|k_z\|_2^{2(s-1)}, \qquad z\in\overline{\dd},
\end{equation}
with the constants depending on $\Theta$ and $s\in (1,\infty)$, but not on
$z$. It follows that 
$$
w_{p,n}(z)\asymp \big(\|k_z\|_2^{2(q(n+1)-1)}\big)^{-\frac{pn}{(pn+1)q}}
=\|k_z\|_2^{-2n}, \quad z\in\dd,
$$
and $w_{p,n}(\zeta) \asymp |\Theta'(\zeta)|^{-n}$, $\zeta\in\tz$.
$\bigcirc$
\medskip
\\
{\bf Remarks.} 1. One should compare inequality (\ref{lev-cls}) 
in Theorem 1.6 with a Bernstein inequality for $L^\infty$-norms:
if $\Theta$ has a nontagential limit at the point
$\zeta\in \mathbb{T}$, $\Theta(\zeta) \in \tz$
and $|\Theta'(\zeta)|<\infty$, then
for each $f\in K_\Theta^\infty$ the derivative 
$f'(\zeta )$ exists in the sense
of non-tangential boundary values and
\begin{equation}
\label{lev}
|f'(\zeta )/\Theta'(\zeta )| \le \|f\|_\infty.
\end{equation}
Indeed, $f'(\zeta) = \int_{\mathbb{T}} 
\overline{k_\zeta^2(\tau)} f(\tau)\, dm(\tau)$, and hence
$|f'(\zeta)| \le \|f\|_\infty \|k_\zeta\|_2^2= \|f\|_\infty 
|\Theta'(\zeta)|$.
Note that (\ref{lev}) holds for arbitrary (not necessarily $CLS$) 
inner functions, and the constant $1$ is sharp.
This inequality is due to M.B. Levin \cite{l1};
for the case of a finite Blaschke product it was later rediscovered 
by a number of authors (see, e.g., \cite{be}).
\smallskip

2. It is shown in \cite[Example 5.2]{bar05} that the 
exponent $\frac{p}{p+1}$ in the definition of the weight $w_{p,1}$ is,
in a certain sense, best possible.
\smallskip

3. Bernstein inequalities for model subspaces $({K^p_{\theta}})_+$
in the upper half-plane $\mathbb{C}^+$ were previously 
studied by K.M. Dyakonov \cite{d1, d2} who showed
that  differentiation is bounded as an operator
from $({K^p_{\theta}})_+$ to $L^p(\mathbb{R})$ with  $1<p\le\infty$,
that is, \begin{equation}
\label{dyak}
\|f'\|_p\le C(p,\theta)\|f\|_p, \qquad f\in(K^p_{\theta})_+,
\end{equation}
if and only if $\theta'\in H^\infty(\mathbb{C}^+)$.
In this case $\theta$ is meromorphic in the whole complex plane 
and the subspace $(K^p_{\theta})_+$ is closely related 
to a certain space of entire functions  
(in particular, for $p=2$, to de Branges spaces; 
see \cite[Proposition 1.1]{bar03a}).
Weighted Bernstein inequalities of the form (\ref{wei}) 
obtained in \cite{bar03, bar05} essentially generalize inequality 
(\ref{dyak}). An advantage of studying weighted estimates
is that the weight may compensate possible growth of elements of
$(K^p_{\theta})_+$ and their derivatives near the boundary.
Note that the Bernstein inequality for the usual $L^p$-norms,
that is, $\|f'\|_p\le C\|f\|_p$, holds for a model subspace
$K_\Theta^p$ in the disc if and only if it is finite-dimensional
and, thus, $\Theta$ is a finite Blaschke product.
Therefore the idea of weighted Bernstein inequalities
(with an ``improving'' weight) is even more natural when one works
with the spaces in the disc.

\bigskip

\begin{center}
{\bf 3. Embedding theorems}
\end{center}

We show that Theorem 2.1 implies an embedding theorem which generalizes
the Volberg--Treil theorem. This result also gives us a condition sufficient 
for the compactness of the embedding.
In what follows we always assume that $\mu$ is a finite
Borel measure in the closed disc $\overline{\dd}$.
\smallskip

By a {\it square with the sidelength $h$} in the unit disc we mean a set
of the form
\begin{equation}
\label{gsq}
S(h_0,\phi_0, h) =\{\rho e^{i\phi}: h_0 -\frac{h}{2\pi}  \le \rho \le  h_0, \ 
\phi_0 \le\phi \le \phi_0+h \}
\end{equation}
for some $h_0\in (0, 1]$, $\phi_0 \in\rl$, and $0< h< 2\pi h_0$.
By $J(S)$ we denote the {\it lower side} of the square $S$, that is, 
$J(S)= \{h_0 e^{i\phi}: \phi_0 \le\phi \le \phi_0+h \}$.

Note that this definition contains as particular cases 
Carleson squares (\ref{csq}) (they correspond to $h_0=1$)
and dyadic squares (\ref{dsq}) introduced in Section 5.
\smallskip
 
Let $\{S_k\}_{k\in\na}$ be a sequence of squares 
in $\overline{\mathbb{D}}$, let
$J_k$ denote the lower side of the square $S_k$, and let $\delta_{J_k}$
be the Lebesgue measure on the arc $J_k$. Assume that $1<r<p$ and that 
the squares $S_k$ satisfy the following two conditions:
\begin{equation}
\label{so}
\sum\limits_k \delta_{J_k}\in{\cal C}
\end{equation}
and
\begin{equation}
\label{sd}
\sup\limits_{k} |J_k|\cdot \|w_r^{-1}\|_{L^q(J_k)}^p<\infty
\end{equation}
where $w_r(z) = w_{r,1}(z) = \|k_z^2\|_{r'}^{-\frac{r}{r+1}}$, $1/r+1/r'=1$,
is the weight from the Bernstein inequality of Section 2.
Condition (\ref{so}) means that the sequence of squares $\{S_k\}$ 
is sufficiently sparse, whereas their size 
is controlled by inequality (\ref{sd}).
\medskip
\\
{\bf Theorem 3.1.}
{\it Let $\{S_k\}$ be a sequence of squares satisfying (\ref{so})
and (\ref{sd}), and let $\mu$ be a Borel measure on 
$\bigcup\limits_k S_k$. Then 

$(i)$ if $\mu(S_k)\le C|J_k|$, then 
$\mu \in {\cal C}_p(\Theta)$;

$(ii)$ if, moreover, $\mu(S_k)=o(|J_k|)$, then
the embedding $K^p_{\Theta}\subset L^p(\mu)$ is compact.}
\medskip

Note that, as in the Volberg--Treil theorem, we consider the measures with
Carleson-type estimate on a special class of sufficiently large squares.
We will see that the squares in Theorem 3.1 may be 
essentially larger (see Proposition 3.4 \smallskip 
below).

In the proof of Theorem 3.1 we will need the following lemmas.
The first of them shows that the norms of reproducing kernels have  
a certain monotonicity along the radii.
\medskip
\\
{\bf Lemma 3.2.} {\it Let $q>1$. Then there exists
$C=C(q)$ such that for any $z=\rho e^{i\phi}$ and
$\tilde z=\tilde \rho e^{i\phi}$ with $0\le \tilde\rho\le \rho$ we have
\begin{equation}
\label{61-1}
\|k_{\tilde z}\|_q \le C(q) \|k_z\|_q.
\end{equation} }
{\bf Proof.} For the case of the upper half-plane the corresponding property 
is established in \cite[Corollary 4.7]{bar05}. The statement
for the disc follows by the same arguments. $\bigcirc$
\medskip

If the sequence $\{S_k\}$ satisfies (\ref{sd}), then  
it follows from (\ref{61-1}) that
\begin{equation}
\label{61-2}
\sup\limits_{k} |J_k|  
\bigg(\int\limits_{S_k\cap \{|z|=\rho\}} w_r^{-q}(z) 
|dz|\bigg)^{p/q} \le C
\end{equation}
for any $\rho\in (0,1]$ (we denote by $|dz|$ the Lebesgue measure
on the corresponding \medskip 
arc).
\\
{\bf Lemma 3.3.} {\it If $J_k\subset \tz$, then (\ref{sd})  
implies that $\int_{J_k} |\Theta'(\tau)|dm(\tau) <\infty$.
In particular, $Int\, J_k \cap \sigma(\Theta)=\emptyset$
(we denote by $Int\, J_k$ the interior of $J_k$ in $\tz$) 
and $\Theta$ is continuous in each of the $($closed$)$ squares $S_k$. }
\smallskip
\\
{\bf Proof.} By (\ref{geomineq}), inequality (\ref{sd}) implies
$\int_{J_k} |\Theta'(\tau)|^q dm(\tau) <\infty$. Hence,
$\int_{J_k} |\Theta'(\tau)| dm(\tau) <\infty$ and we conclude that
$\Theta$ is continuous on $J_k$. It is easy to see that
$|\Theta'(r\zeta)|\le |\Theta'(\zeta)|$, $\zeta \in \mathbb{T}$,
$r\in (0,1)$, and therefore, $\Theta$ is continuous on $S_k$.
$\bigcirc$
\medskip
\\
{\bf Proof of Theorem 3.1, ({\it i}).}
Clearly, the embedding ${K^p_{\Theta}}\subset L^p(\mu|_{\{|z|< 1/2\}})$
is compact. So we may assume without loss of generality 
that $supp \, \mu\subset \{1/2 \le |z|\le 1\}$.
It follows from Lemma 3.3 that the set of functions $f\in{K^p_{\Theta}}$
which are continuous on each of $S_k$
is dense in ${K^p_{\Theta}}$, $1<p<\infty$ (take the reproducing kernels). 
Thus it is sufficient to prove the estimate 
$$
\|f\|_{L^p(\mu)} \le C\|f\|_p, \qquad f\in{K^p_{\Theta}},
$$
only for $f$ continuous in $\bigcup\limits_k S_k$.

Now let $f \in K^p_{\Theta}$ be continuous in each of $S_k$. 
Then there exist $w_k \in S_k$
such that
\begin{equation}
\label{61-4}
\|f\|^p_{L^p(\mu)} \le \sum\limits_k |f(w_k)|^p\mu(S_k) \le
\sup\limits_k \frac{\mu(S_k)}{|J_k|} \cdot \sum\limits_k |f(w_k)|^p |J_k|.
\end{equation}
Statement ($i$) will be proved as soon as we show that
\begin{equation}
\label{61-5}
\sum\limits_k |f(w_k)|^p |J_k| \le C\|f\|_p^p
\end{equation}
where $C$ does not depend on $f$ and on the choice of $w_k\in S_k$.

Consider the arcs $\tilde J_k = S_k\cap \{|z|= |w_k|\}$. Since 
$\mu(\{|z|<1/2\})=0$, it follows that $|\tilde J_k| \ge |J_k|/2$.
Let $\nu =\sum_k \delta_{\tilde J_k}$. Then it follows from (\ref{so})
that $\nu\in{\cal C}$ (and the Carleson constants $M_\nu$ of 
such measures $\nu$
are uniformly bounded). We have
\begin{equation}
\label{61-3}
\bigg(\sum\limits_k |f(w_k)|^p |\tilde J_k|\bigg)^{1/p}
\le \|f\|_{L^p(\nu)} +
\bigg(\sum\limits_k \int\limits_{\tilde J_k} |f(z)-f(w_k)|^p |dz|\bigg)^{1/p},
\end{equation}
and $\|f\|_{L^p(\nu)} \le C_1 \|f\|_p$. 

We estimate the last term in (\ref{61-3}).
For $z\in \tilde J_k$ 
denote by $\gamma(z, w_k)$ the subarc of $\tilde J_k$ with the endpoints
$z$ and $w_k$.
Then $f(z)-f(w_k) = \int_{\gamma(z, w_k)} f'(u) du$
(in the case $J_k \subset \tz$ note that, by Lemma 3.3,  
any $f\in {K^p_{\Theta}}$ is analytic on $J_k$ 
except, may be, at the endpoints), and so
$$
\sum\limits_k \int\limits_{\tilde J_k} |f(z)-f(w_k)|^p |dz| =
\sum\limits_k \int\limits_{\tilde J_k} \bigg| 
\int\limits_{\gamma(z, w_k)} f'(u) du\bigg|^p |dz|
$$ 
$$
\le \sum\limits_k \int\limits_{\tilde J_k} 
\bigg(\int\limits_{\gamma(z, w_k)} w_r^{-q}(u) |du|\bigg)^{p/q}
\bigg(\int\limits_{\gamma(z, w_k)} |f'(u)|^p w_r^p(u) |du|\bigg) |dz|
$$
$$
\le  \sum\limits_k |\tilde J_k| 
\bigg(\int\limits_{\tilde J_k} w_r^{-q}(u) |du|\bigg)^{p/q}
\bigg(\int\limits_{\tilde J_k} |f'(u)|^p w_r^p(u) |du|\bigg).
$$
By (\ref{61-2}), we obtain
$$
\sum\limits_k \int\limits_{\tilde J_k} |f(z)-f(w_k)|^p |dz|
\le C_2 \sum\limits_k \int\limits_{\tilde J_k} |f'(u)|^p w_r^p(u) |du|
$$
$$
= C_2 \|f' w_r\|^p_{L^p(\nu)} \le C_3\|f\|_p^p
$$
where the last inequality follows from Theorem 2.1. 
\smallskip
\\
{\bf Proof of Statement ({\it ii}).} For a Borel set
$E\subset\overline{\mathbb{D}}$ define the operator
${\cal I}_E: {K^p_{\Theta}}\to L^p(\mu)$ by ${\cal I}f= \chi_E f$
where $\chi_E$ is the characteristic function of $E$.
For $N\in \na$ put $F_N= \bigcup\limits_{k=1}^N S_k$. As above we assume that
$f\in {K^p_{\Theta}}$ is continuous in $\overline{\mathbb{D}}$. 
Then it follows from (\ref{61-4}) and (\ref{61-5}) that  
$$
\int\limits_{\overline{\mathbb{D}}\setminus F_N} |f|^p d\mu \le C\sup\limits_{k>N}
\frac{\mu(S_k)}{|J_k|}\|f\|_p^p,
$$
and so $\|{\cal I}_{\overline{\mathbb{D}}\setminus F_N}\| \to 0$, 
$N\to\infty$. Statement ($ii$) will be proved as soon as we show that
${\cal I}_{F_N}$ is a compact operator for any $N$ (thus,
our embedding operator ${\cal J}_\mu= {\cal I}_{\overline{\mathbb{D}}}=
{\cal I}_{F_N}+{\cal I}_{\overline{\mathbb{D}}\setminus F_N}$
may be approximated in the operator norm by compact operators 
${\cal I}_{F_N}$). Clearly, it suffices to prove the compactness 
of ${\cal I}_{S_k}$ for each fixed $k$.

We approximate ${\cal I}_{S_k}$ by finite rank operators.
Partition the square $S_k$ into finite union of 
pairwise disjoint squares $\{\tilde S_l\}_{l=1}^L$
(here we do not require $\tilde S_l$ to be closed and assume that
$Clos\, \tilde S_l$ is a closed square in the sense 
of definition (\ref{gsq})) and choose a point $\zeta_l$ 
in each of $\tilde S_l$. By Lemma 3.2, we can choose 
the squares $\tilde S_l$ to be so small that, for a given $\epsilon>0$,
we would have
\begin{equation}
\label{61-8}
\bigg(\, \int \limits_{[\zeta_l, z_l]} 
w^{-q}_r(z) |dz| \bigg)^{p/q} <\epsilon
\end{equation}
for any $l$, $1\le l\le L$, and any $z_l\in \tilde S_l$.
Here we denote by $[z,w]$ the straight line interval with 
endpoints $z$ and $w$.

Now we consider the finite rank operator $T: {K^p_{\Theta}} \to L^p(\mu)$,
$(Tf)(z) = \sum_{l=1}^L f(\zeta_l) \chi_{\tilde S_l}(z)$.
We show that $\| {\cal I}_{S_k} - T\|^p \le C\epsilon$. As in the proof
of ($i$), we have
$$
\|({\cal I}_{S_K} - T) f\|_{L^p(\mu)}^p = 
\sum\limits_{l=1}^L \int\limits_{\tilde S_l}|f(z)-f(\zeta_l)|^pd\mu(z)
$$
$$
\le 
\sum\limits_{l=1}^L   \int\limits_{\tilde S_l}
\bigg(\, \int \limits_{[\zeta_l, z]} |f'(u)|^p w_r^p(u) |du| \bigg)\cdot
\bigg(\, \int \limits_{[\zeta_l, z]} w^{-q}_r(u) |du| 
\bigg)^{p/q} d\mu(z).
$$
By Theorem 2.1,
$$
 \int \limits_{[\zeta_l, z]} |f'(u)|^p w^p_r(u) |du| \le C_1 \|f\|_p^p
$$
where $C_1$ does not depend on $f\in{K^p_{\Theta}}$, $1\le l\le L$ and $z\in \tilde S_l$.
Hence, by (\ref{61-8}), 
$$
\|({\cal I}_{S_K} - T) f \|_{L^p(\mu)}^p \le C_1\epsilon 
\|f\|_p^p
\sum\limits_{l=1}^L \mu(\tilde S_l) = C_1\epsilon \mu(S_k) \|f\|_p^p.
$$
We conclude that ${\cal I}_{S_K}$ may be approximated by finite rank operators
and therefore is compact. $\bigcirc$
\medskip
\\
{\bf Remark.} In the proof of Theorem 3.1 we used the fact that,
by Lemma 3.3, the functions continuous in each of $S_k$
are dense in $K^p_\Theta$. Actually, it follows from the results 
of \cite{al9} that the functions continuous in the closed disc 
$\overline{\mathbb{D}}$ are dense in $K^p_\Theta$, $p\ge 1$.
\medskip

Theorem 3.1 describes a wider class of measures than the 
class ${\cal C}(\Theta)$ in the Volberg--Treil theorem.
Namely, we have the following proposition:
\medskip
\\
{\bf Proposition 3.4.}
{\it If $\mu\in{\cal C}(\Theta)$, then $\mu=\mu_1+\mu_2$ where 
$\mu_1$ satisfies the conditions of Theorem 3.1 $(i)$ for any $p>1$ 
and $r\in (1,p)$, whereas $\mu_2\in{\cal C}$. }
\medskip

In what follows we will use a special family of arcs on $\tz$
(analogous to the Whitney decomposition for the set 
$\overline{\dd}\setminus \Omega(\Theta,\varepsilon)$).
\medskip
\\
{\bf Lemma  3.5.} {\it
Let $\varepsilon\in(0,1)$.
Assume that $\tz\setminus \sigma(\Theta) \ne \emptyset$. Then there exists
a sequence of arcs $I_k\subset \tz$, $k\in \na$,
with pairwise disjoint interiors 
such that $\bigcup\limits_k I_k =
\tz\setminus \sigma(\Theta) $ and
\begin{equation}
\label{tgv}
|I_k| \le dist (I_k, \Omega(\Theta,\varepsilon)) \le 2|I_k|.
\end{equation}
Moreover, if we put $F=\bigcup \limits_k S(I_k)$ and $G= \overline{\dd} 
\setminus F$, then for any $z\in G$, $z\ne 0$,
we have 
\begin{equation}
\label{tgv1}
dist (z/|z|,\Omega(\Theta,\varepsilon)) \le 6\pi (1-|z|).
\end{equation}}
{\bf Proof.} Note that 
$\int_{\tz\setminus \sigma(\Theta)} d_\varepsilon^{-1}(\zeta)
\, dm(\zeta) =\infty$.
Therefore we can choose the sequence of arcs $I_k$ 
with pairwise disjoint interiors 
such that
$\bigcup\limits_k I_k =
\tz\setminus \sigma(\Theta) $
and 
$$
\int\limits_{I_k} d_\varepsilon^{-1}(\zeta)\, dm(\zeta) =\frac{1}{2}.
$$
It follows that there exists $\zeta_k\in I_k$ such that
$d_\varepsilon (\zeta_k) =2|I_k|$. Hence, for any $\zeta\in I_k$,
$d_\varepsilon (\zeta) \ge d_\varepsilon(\zeta_k) - |I_k| \ge |I_k|$,
and we get (\ref{tgv}).

Now let $z=re^{i\phi} \in G$. Then, either $e^{i\phi} \in \sigma(\Theta)$ 
(and so $dist (e^{i\phi},\Omega(\Theta,\varepsilon)) \le 1-|z|$), 
or $e^{i\phi} \in I_k$ for some $k$.
Since $z\notin S(I_k)$, we have $1-r\ge |I_k|/(2\pi)$. Hence, by (\ref{tgv}),
$dist (e^{i\phi}, \Omega(\Theta,\varepsilon))\le 3|I_k| \le 6\pi (1-r)$. $\bigcirc$
\medskip
\\
{\bf Proof of Proposition 3.4.} 
As usual, for an arc $I\subset \tz$ and $a>0$ we denote by $aI$
the arc with the same center of the length $a|I|$.
Put $\mu_1=\mu|_F$ and  $\mu_2=\mu|_G$ where the sets 
$F, G$ are defined in Lemma 3.5.
It follows from (\ref{tgv}) that
$$
|I_k| \bigg(\int\limits_{I_k} 
d_\varepsilon^{-q}(\zeta)\, dm(\zeta)\bigg)^{p/q} \le C.
$$
Let $r\in (1,p)$. By (\ref{geomineq}), 
$d_\varepsilon\lesssim w_r$, and therefore we have 
\begin{equation}
\label{tgv5}
|I_k| \bigg(\int\limits_{I_k} 
w_r^{-q}(\zeta)\, dm(\zeta)\bigg)^{p/q} \le C_1.
\end{equation}
Thus, the family of the squares $S(I_k)$ satisfies (\ref{so}) and (\ref{sd}).
Let us show that $\mu_1(S(I_k)) \le C_2|I_k|$. Indeed, 
it follows from (\ref{tgv})
that $S(A I_k)\cap \Omega(\Theta,\varepsilon)\ne \emptyset$ for some sufficiently large 
absolute constant $A>1$. Since $\mu\in {\cal C}(\Theta)$, it follows
that $\mu_1 (S(I_k)) \le \mu (S(AI_k)) \le C_2 |I_k|$.

Now we show that $\mu_2$ is a usual Carleson measure. 
If $S(I) \cap G\ne\emptyset$
for some arc $I\subset \tz$, then, by (\ref{tgv1}),
there is an absolute constant $A_1>1$ such that
\begin{equation}
\label{l53}
S(A_1 I) \cap \Omega(\Theta,\varepsilon)\ne\emptyset,
\end{equation}
and so $\mu(S(I)) \le C_3|I|$ for a positive constant $C_3$.
$\bigcirc$
\medskip

The following example shows that Theorem 3.1 describes  
an essentially wider class
of embeddings than the Volberg--Treil theorem.
\medskip
\\
{\bf Example 3.6.} By Proposition 3.4,
each measure $\mu\in{\cal C}(\Theta)$ is of the form $\mu=\mu_1+\mu_2$ where $\mu_1$
satisfies condition ($i$) of Theorem 3.1 and $\mu_2$ is a usual Carleson measure.
Thus, the Volberg--Treil theorem follows from Theorem 3.1.
On the other hand, it is easy to construct
a measure $\mu$ satisfying the conditions of
Theorem 3.1 ($i$), which is not in ${\cal C}(\Theta)$.

Clearly, if $\mu\in{\cal C}(\Theta)$, then $\mu$ has no point masses at the points of
the boundary spectrum, that is, $\mu(\{\zeta\})=0$ for any $\zeta\in\sigma(\Theta)
\cap \tz$ (note that in this case $S(I)\cap\Omega(\Theta,\varepsilon)
\ne\emptyset$ for any arc $I$ such that $\zeta$ is an interior point of $I$
and for any $\varepsilon\in (0,1)$).
However, measures in Theorem 3.1 can have nontrivial point masses on
$\sigma(\Theta) \cap \tz$. Let $B$ be the Blaschke product with the zeros
$z_n= (1-2^{-n})e^{i/n}$, $n\in\na$. Then for any $p\in(1,\infty)$
we have 
$$
\|k_\zeta^2\|_q\le C_1, \qquad -\pi <\arg\zeta\le 0,
$$  
and, consequently, $w_p^{-1}(\zeta)=\|k_\zeta^2\|_q^{p/(p+1)}\le C_2$,
$-\pi <\arg\zeta\le 0$.
Hence, $\delta_1\in {\cal C}_p(\Theta)$. Analogously, it is easy to construct
an infinite sum of point masses, that is
$\mu=\sum_n a_n \delta_{\zeta_n}$ with $a_n>0$,
$\zeta_n\in\sigma(\Theta) \cap \tz$, such that the embedding
${K^p_{\Theta}}\subset L^p(\mu)$ is bounded or even compact 
(see \cite[Example 6.3]{bar05} for details).

\bigskip

\begin{center}
{\bf 4. Compact embeddings. Proof of Theorem 1.1}
\end{center}

In this section we prove Theorem 1.1 and discuss the relation between 
the two ``vanishing conditions'' introduced in \cite{cimmat}.
\medskip
\\
{\bf Proof of Theorem 1.1.}
As we mentioned in Introduction, implication 
($ii$)$\Longrightarrow$($i$) for $CLS$ inner functions 
was proved in \cite{cimmat}. We show that
($i$)$\Longrightarrow$($ii$); the proof is analogous 
to the proof of Proposition 3.4.
Let $\varepsilon\in (0,1)$ and let $I_k$, $F$, $G$, 
$\mu_1$ and $\mu_2$ have the same meaning as in 
Proposition 3.4 and Lemma 3.5. We show that
$\mu_1$ satisfies condition ($ii$) of Theorem 3.1 and $\mu_2$
is a vanishing Carleson measure (see Section 2 for the definition),
and thus the embedding $H^p\subset L^p(\mu_2)$ is compact.

For any $p\in (1,\infty)$ and $r\in(1,p)$, the inequality 
(\ref{tgv5}) is satisfied. Since $|I_k| \to 0$, $k\to\infty$, 
and $S(AI_k)\cap \Omega(\Theta,\varepsilon)\ne \emptyset$
for a sufficiently large 
absolute constant $A>1$, it follows from the hypothesis that
$$
\lim\limits_{k\to \infty}\frac{\mu(S(I_k))}{|I_k|}=0.
$$
Hence the embedding $K_\Theta^p \subset L^p(\mu_1)$
is compact by Theorem 3.1, ($ii$). As we have shown 
in the proof of Proposition 3.4, for any arc $I$ such that $\mu_2(S(I))\ne 0$
(that is, $S(I) \cap G \ne \emptyset$),
we have (\ref{l53}) for a sufficiently large absolute constant $A_1>1$. 
By condition ($i$) of Theorem 1.1, $\mu(S(I))/|I|\to 0$ when $|I|\to 0$ and
$S(A_1I)\cap \Omega(\Theta,\varepsilon)\ne\emptyset$. 
Hence $\mu_2$ is a vanishing Carleson measure. $\bigcirc$
\medskip

In \cite{cimmat} another sufficient  condition 
for the compactness of the embedding 
was introduced. Following \cite{cimmat}, we put, for $\delta>0$,
$$
H_\delta=\{z\in\overline{\dd}: dist (z, \sigma(\Theta)\cap \tz) <\delta\}
$$
and we say that a measure $\mu$ satisfies the {\it first vanishing condition} 
(V1, for short) if 
\begin{equation}
\label{71-2}
M_{\mu_\delta} \to 0, \qquad \delta\to 0,
\end{equation}
where $\mu_\delta=\mu|_{H_\delta}$. Recall that $M_\nu$ denotes
the Carleson constant of a  Carleson measure $\nu$. 
If a measure $\mu$ satisfies condition ($i$) of Theorem 1.1, 
we say that $\mu$ satisfies the {\it second vanishing condition} (V2).

It is shown in  \cite{cimmat} that if $\mu$ satisfies V1, then the embedding
${K^p_{\Theta}} \subset L^p(\mu)$ is compact for  $0<p<\infty$.
The authors ask what the relation between 
the two vanishing conditions is. Here we 
answer this question by showing that V1 always implies V2, but
not vice versa (thus, for $p > 1$, the sufficient condition of compactness
given in \cite{cimmat} follows from our Theorem 1.1).
\medskip
\\
{\bf Proposition 4.1.} {\it V1 implies V2.}
\smallskip
\\
{\bf Proof.}
Assume that $\mu$ satisfies V1, but not V2. 
Then there exists a sequence of arcs $\{J_n\}_{n\in \na}$
such that, for some fixed $\varepsilon\in(0,1)$, $S(J_n)\cap 
\Omega(\Theta,\varepsilon) \ne \emptyset$, $|J_n|\to 0$, $n\to\infty$, but 
\begin{equation}
\label{71-0}
\mu(S(J_n))\ge C|J_n|
\end{equation}
for a constant $C>0$.

Fix $\delta>0$ and put $G_\delta =\overline{\dd}\setminus H_\delta$.
It follows from the definition of $\sigma(\Theta)$ that $\Theta$ 
is continuous on $G_\delta$ and 
$|\Theta(z)|\to 1$ uniformly when $|z|\to 1$, $z\in G_\delta$. 
Hence there exists a constant $\delta_1 \in (0,\delta)$ such that
\begin{equation}
\label{71-1}
|\Theta(z)| > \varepsilon, \qquad z\in G_\delta, \ \  1-\delta_1 \le |z| \le 1.
\end{equation}

Choose $N$ such that  $|J_n| <\delta_1$, $n \ge N$. 
Clearly, $S(J_n) \subset \{z \in \overline{\dd}: \, 1-\delta_1 \le |z| \le 1 \}$.
Since $S(J_n) \cap \Omega(\Theta,\varepsilon) \ne \emptyset$, it follows from (\ref{71-1}) that 
$$
S(J_n)  \nsubseteq G_\delta \cap \{1-\delta_1 \le |z|\le 1\}, \qquad n\ge N.
$$ 
We conclude that
$S(J_n)\cap H_\delta \ne \emptyset$, $n\ge N$.
Hence $S(J_n) \subset H_{2\delta}$, and, consequently, 
$\mu(S(J_n)) \le M_{\mu_{2\delta}} |J_n|$, $n\ge N$. 
This contradicts (\ref{71-0}), since, by (\ref{71-2}), 
$ M_{\mu_{2\delta}} \to 0$ when $\delta\to 0$. $\bigcirc$
\medskip
\\
{\bf Example 4.2.}
Now we show that V1 is not necessary for the compactness of the embedding
even for $CLS$ inner functions (and, thus, V2 does not imply V1).
Let $I_n$ be the sequence of arcs from Lemma 3.5, and let
$\zeta_n$ be the middle point of $I_n$. 
Put $\mu=\sum_n a_n |I_n|\delta_{\zeta_n}$
where $a_n\to 0$, $n\to\infty$. By Theorem 3.1, ($ii$), 
the embedding ${K^p_{\Theta}} \subset L^p(\mu)$ is compact for any 
$p\in(1,\infty)$. 
Assume that $\Theta \in CLS$. Now, by Theorem 1.1 ($ii$),
$\mu$  satisfies V2 (this can be also 
shown directly by the arguments analogous
to those in the proof of Proposition 3.4). However, the measure $\mu_\delta$
has nonzero point masses on $\mathbb{T}$ for any $\delta>0$, and hence 
is not a Carleson measure. Thus, $\mu$ does not satisfy V1.

\bigskip

\begin{center}
{\bf 5. Classes ${\cal S}_r$. Sufficient conditions}
\end{center}

The definition and basic properties of Schatten--von Neumann
operator ideals ${\cal S}_r$ may be found in \cite{gohb}.
\smallskip

By $R_{n,m}$ we denote the elements of the 
standard dyadic partition of $\dd$; namely,
\begin{equation}
\label{dsq}
R_{n,m} =\left\{z= \rho e^{i\phi}: 
1-\frac{1}{2^{n-1}} \le \rho < 1-\frac{1}{2^n}, \
\frac{\pi m}{2^{n-1}} \le \phi <\frac{\pi (m+1)}{2^{n-1}} \right\}
\end{equation}
where $n \in \na$, $m=0, 1,\dots 2^n-1$.

We recall a theorem due to Luecking \cite{luek}
concerning Schatten--von Neumann
properties of the embeddings of the whole Hardy space $H^2$. 
For a measure $\mu \in{\cal C}$ 
the embedding operator of $H^2$ into $L^2(\mu)$ 
is in ${\cal S}_r$, $r>0$, if and only if 
\begin{equation}
\label{823-0}
\sum\limits_{n,m} (2^{n}\mu(R_{n,m}))^{r/2}<\infty
\end{equation}
where we sum over all dyadic squares $R_{n,m}$.
An interesting general approach to embeddings
of reproducing kernel Hilbert spaces was suggested 
by Parfenov \cite{parf1}. We will essentially use
the ideas from \cite{parf1}, especially in the proof of necessary 
\smallskip
conditions.  

A criterion for the inclusion of the embedding operator ${\cal J}_\mu:
{K^2_\Theta} \to L^2(\mu)$, ${\cal J}_\mu f=f$, into ${\cal S}_2$ is obvious.
\medskip
\\
{\bf Proposition 5.1.} ${\cal J}_\mu \in {\cal S}_2$ if and only if
$\|k_z\|_2 \in L^2(\mu)$; moreover,
$\|{\cal J}_\mu \|^2_{{\cal S}_2} =\int \|k_z\|_2^2 \, d\mu(z)$.
\smallskip
\\
{\bf Proof.} We have
\begin{equation}
\label{proekt}
({\cal J}_\mu f)(z) =
\int\limits_\tz f(w)\overline{k_z(w)} \,dm (w),\qquad f\in{K^2_\Theta}.
\end{equation}
Note that the operator $\widetilde{{\cal J}}_\mu$ defined by (\ref{proekt}) 
on the whole $L^2(\tz)$ is the orthogonal projection
from $L^2(\tz)$ onto ${K^2_\Theta}$. Then ${\cal J}_\mu \in {\cal S}_2$
if and only if $\widetilde{{\cal J}}_\mu \in {\cal S}_2$
which is equivalent to
$$
\int \int\limits_{\tz} |k_z(w)|^2 \, dm(w)\, d\mu(z) =
\int \|k_z\|_2^2 \, d\mu(z)<\infty. \qquad \bigcirc
$$

Take the point $z_{n,m}=(1-2^{-n})
\exp (i\pi 2^{-n} m)$ in $R_{n,m}$. 
It is easy to see that $\|k_z\|_2^2 \asymp 2^{n}(1-|\Theta(z_{n,m})|)$,
$z\in R_{n,m}$, with the constants independent of $n,m$. Then
the condition $\|k_z\|_2 \in L^2(\mu)$ may be rewriten as
$$
\sum\limits_{n,m} 2^n \mu(R_{n,m}) (1-|\Theta(z_{n,m})|) <\infty.
$$

Our next result in this section is a 
sufficient condition for the inclusion ${\cal J}_\mu \in {\cal S}_r$, $r>0$,
which involves the arcs $I_k$ from Lemma 3.5 and 
some special families of dyadic squares; it contains Theorem 1.2.
For $\varepsilon \in (0,1)$ and $A>0$ put 
$$
{\cal R}(\varepsilon, A)= \{R_{n,m}: dist (R_{n,m},\Omega(\Theta,\varepsilon)) \le A2^{-n}\}.
$$
{\bf Theorem 5.2.}
{\it Let $r>0$, let $\mu$ be a Borel measure in $\overline{\dd}$,
and let $\varepsilon \in (0,1)$. There exists
an absolute constant $A>0$ such that 
${\cal J}_\mu \in {\cal S}_r$ whenever
\begin{equation}
\label{23-00}
\sum\limits_{k} \left(\frac{\mu(S(I_k))}{|I_k|}\right)^{r/2}<\infty
\end{equation}
and
\begin{equation}
\label{m1.8}
\sum\limits_{R_{n,m}\in {\cal R}(\varepsilon, A)} 
(2^n \mu(R_{n,m}))^{r/2}<\infty.
\end{equation} }
\smallskip

We introduce the following operator $T_\mu: K_\Theta^2 \to K_\Theta^2$,
\begin{equation}
\label{tmu}
(T_\mu f)(w) =\int f(z) \overline{k_w(z)}\, d\mu (z).
\end{equation}

The following lemma follows from (\ref{proekt}) 
by a simple calculation:
\medskip
\\
{\bf Lemma 5.3.}
{\it We have $ {\cal J}_\mu^*{\cal J}_\mu = T_\mu $ 
and for any $f,g \in K_\Theta^2$,
$$
\langle T_\mu f, g \rangle  =\int f(z)\overline{g(z)}\, d\mu(z)
$$
where $\langle\cdot, \cdot\rangle$ stands for the 
standard inner product in $L^2(\tz)$.}
\medskip

In the proof of Theorem 5.2 we will use the following property
of the arcs $I_n$ constructed in Lemma 3.5.
\medskip
\\
{\bf Lemma 5.4.} {\it Let $\varepsilon\in (0,1)$, 
and let $\{I_n\}$ be the system of arcs
of Lemma 3.5. Then there exists a constant $C=C(\varepsilon)>0$
such that
$$
|k_z(w)| = \left|\frac{1-\overline{\Theta(z)}\Theta(w)}
{1-\overline z w} \right| \le C|I_n|^{-1}
$$
for any $n$ and $z,\, w \in S(I_n)$. In particular, $|\Theta'(\zeta)|\le 
C|I_n|^{-1}$, $\zeta\in I_n$.}
\medskip
\\
{\bf Proof.} By construction of $I_n$,
$d_\varepsilon(\zeta) = dist\, (\zeta, \Omega(\Theta,\varepsilon)) \ge |I_n|$,
$\zeta \in I_n$. 
Let $w=r\zeta$, $r\in (0,1)$, $\zeta\in I_n$. Clearly,
$|k_z(w)| \le \|k_z\|_2\|k_w\|_2$, 
and, by Lemma 3.2, 
$$
\|k_w\|_2^2\le C \|k_\zeta\|_2^2=C|\Theta'(\zeta)|.
$$
It follows from (\ref{geomineq}) (see also \cite[Theorem 4.9]{bar05}) 
that $|\Theta'(\zeta)| \le C_1 (d_\varepsilon(\zeta))^{-1}
\le C_1 |I_n|^{-1}$. Hence, $\|k_w\|_2^2 \le C_2 |I_n|^{-1}$, 
$w\in S(I_n)$.  $\bigcirc$
\medskip

Now we turn to the proof of Theorem 5.2. We start 
with the proof for the case $r\ge 2$;
here we follow the argument from \cite{luek} based on the  method of 
complex interpolation between the different Schatten--von Neumann 
classes. Then we use an idea from
Parfenov's paper \cite{parf1} to give the proof for 
$0<r\le 1$. Finally, the case $1<p<2$ follows by interpolation.
\medskip
\\
{\bf Proof of Theorem 5.2.} Given $\varepsilon \in (0,1)$, consider 
the system of arcs $\{I_n\}$ constructed in Lemma 3.5.
As in the proof of Proposition 3.4, 
we put $\mu_1=\mu|_F$ and  $\mu_2=\mu|_G$ where
$F=\bigcup \limits_n S(I_n)$ and $G= \overline{\dd}  \setminus F$.

First we show that, for any $r>0$,
condition (\ref{m1.8}) with appropriate $A$
implies that
the embedding operator ${\cal J}_2: H^2\to L^2(\mu_2)$ is in ${\cal S}_r$.
Let $R_{n,m}$ be a dyadic square
such that $R_{n,m} \cap G\ne \emptyset$.
We show that $R_{n,m}\in {\cal R}(\varepsilon, A)$ for some $A>0$.
Indeed, let $z\in R_{n,m} \setminus F$, 
$z=(1-\rho) \zeta$, $\rho\in(0,1)$, $\zeta\in\tz$.
If $\zeta\in \sigma(\Theta)$, then 
$dist (z,\Omega(\Theta,\varepsilon)) \le \rho \le 2^{-(n-1)}$.
Otherwise, $\zeta\in I_k$ and, by
(\ref{tgv1}), $dist (\zeta,\Omega(\Theta,\varepsilon)) \le 6\pi\rho$.
Hence, $dist (z,\Omega(\Theta,\varepsilon)) \le (6\pi+1)\rho \le (6\pi+1) 2^{-(n-1)}$.
We conclude that $R_{n,m}\in {\cal R}(\varepsilon, A)$ with $A=12\pi+2$.
Hence
$$
\sum\limits_{R_{n,m}\cap G \ne \emptyset} (2^{n}\mu(R_{n,m}))^{r/2}<\infty,
$$
and ${\cal J}_2 \in {\cal S}_r$ by the Luecking theorem.
\smallskip

Now we consider the embedding operator 
${\cal J}_1: {K^2_\Theta} \to L^2(\mu_1)$.
\smallskip
\\
{\bf Proof for the case $r\ge 2$.} 
By Lemma 5.3, ${\cal J}_1 \in {\cal S}_r$ if and only if
the operator
$$
(Tf)(w) = ({\cal J}_1^* {\cal J}_1 f)(w) = 
\int\limits_{F} f(z) \overline{k_w(z)}\, d\mu (z)
$$
is in ${\cal S}_{r/2}$. Let $p=r/2$. Since $r\ge 2$ we have $p\ge 1$. 
For $\zeta\in \co$, $0\le {\rm Re}\, \zeta \le 1$, we put
$$
(T(\zeta)f) (w) =
\sum\limits_n \left(\frac{\mu(S_n)}{|I_n|}\right)^{\zeta p-1}
\int\limits_{S_n} f(z) \overline{k_w(z)} d\mu (z)
$$
where $S_n =S(I_n)$. Then $T(\zeta)$ is an analytic family 
of operators in ${K^2_\Theta}$.
Clearly, $T(1/p)=T$. We will show that $T(\zeta)$ is bounded in ${K^2_\Theta}$
with $\|T(\zeta)\|\le A_0$, ${\rm Re}\, \zeta\in [0,1]$.
Also we will show that
\begin{equation}
\label{83-30}
\|T(\zeta)\|_{{\cal S}_1} \le A_1, \qquad {\rm Re}\, \zeta =1,
\end{equation}
where $\|T(\zeta)\|_{{\cal S}_1}$ is the trace norm of the operator $T(\zeta)$.
Then, by \cite[Theorem 13.1]{gohb}, 
$T=T(1/p) \in {\cal S}_p$ and
$\|T\|_{{\cal S}_p} \le A_0^{1-1/p}A_1^{1/p}$.

Note that for any $f,g\in{K^2_\Theta}$ we have, as in Lemma 5.3,
\begin{equation}
\label{83-21}
\langle T(\zeta)f, g\rangle = \sum\limits_n
\left(\frac{\mu(S_n)}{|I_n|}\right)^{\zeta p-1}
\int\limits_{S_n}  f(z) \overline{g(z)}\,d\mu(z).
\end{equation}
Note also that, by (\ref{23-00}), $\mu(S_n) \le C|I_n|$, and so 
$\big|(\mu(S_n)/|I_n|)^{\zeta p-1}\big| \le 
C_1 (\mu(S_n)/|I_n|)^{-1}$ for any $\zeta$
with ${\rm Re}\, \zeta\in [0,1]$.
Hence, 
$$
|\langle T(\zeta)f, g\rangle |\le C_1 \sum\limits_n
(\mu(S_n))^{-1}|I_n|\int\limits_{S_n}
|f(z)g(z)|d\mu(z)
$$
$$
\le \sum\limits_n |I_n|\cdot |f(z_n)g(z_n)|\le
\bigg(\sum\limits_n |I_n|\cdot |f(z_n)|^2\bigg)^{1/2} 
\bigg(\sum\limits_n |I_n|\cdot |g(z_n)|^2\bigg)^{1/2} 
$$
for some points $z_n\in S_n$ (recall that, by Lemma 3.3, 
the functions continuous in $S_n$ are dense in
$K^2_\Theta$). By Theorem 3.1, 
the measure  $\nu=\sum\limits_n |I_n|\,\delta_{z_n}$
is in the class ${\cal C}_2(\Theta)$
(indeed, $\nu(S_n) = |I_n|$), 
and, moreover, Carleson constants of such measures
are uniformly bounded by a constant $C_2$ which does not depend
on the choice of $z_n\in S_n$. Thus,
$$
\sum\limits_n |I_n|\cdot |f(z_n)|^2\le C_2^2\|f\|_2^2, \qquad f\in{K^2_\Theta}.
$$
We conclude that
$$
|\langle T(\zeta)f, g\rangle |\le C_3 \|f\|_2  \|g\|_2,\qquad f,g\in{K^2_\Theta},
$$
which implies $\|T(\zeta)\|\le A_0$, ${\rm Re}\, \zeta\in [0,1]$.
\smallskip

It remains to verify (\ref{83-30}).
Note that condition (\ref{23-00}) implies
$\mu(S_n)=o(|I_n|)$, $n\to\infty$, and therefore $T(\zeta)$
is compact for ${\rm Re}\, \zeta =1$. Hence, 
for a fixed $\zeta$, we can write $T(\zeta)$
in the canonical form $T(\zeta) h=\sum_m  \beta_m
\langle h, e_m\rangle f_m$
where $\{e_m\}$ and $\{f_m\}$ are orthonormal systems in
${K^2_\Theta}$,
$\{\beta_m\}$ is a sequence of positive numbers tending to zero,
and $\|T(\zeta)\|_{{\cal S}_1} =\sum_m \langle T(\zeta) e_m, f_m\rangle$.
By (\ref{83-21}), 
$$
\sum\limits_m \langle T(\zeta)e_m, f_m\rangle=
\sum\limits_m \bigg(\sum\limits_n
\left(\frac{\mu(S_n)}{|I_n|}\right)^{\zeta p-1}
\int\limits_{S_n}  e_m(z) \overline{f_m(z)}\,d\mu(z)\bigg)
$$
$$
\le \sum\limits_n
\left(\frac{\mu(S_n)}{|I_n|}\right)^{p-1}
\bigg(\int\limits_{S_n} \sum\limits_m |e_m(z)|^2 d\mu(z)\bigg)^{1/2}
\bigg(\int\limits_{S_n} \sum\limits_m |f_m(z)|^2 d\mu(z)\bigg)^{1/2}.
$$
By Parseval's identity, $\sum_m |g_m(z)|^2 \le k_z(z)=\|k_z\|_2^2$
for any orthonormal system $\{g_m\}$ in ${K^2_\Theta}$. Thus, 
for ${\rm Re}\, \zeta=1$,
$$
\|T(\zeta)\|_{{\cal S}_1} \le \sum\limits_n
\left(\frac{\mu(S_n)}{|I_n|}\right)^{p-1}
\int\limits_{S_n} \|k_z\|_2^2 \,  d\mu(z).
$$
By Lemma 5.4, $\|k_z\|_2^2=k_z(z) \le C|I_n|^{-1}$, $z\in S_n$. 
Hence
$$
\|T(\zeta)\|_{{\cal S}_1} \le C \sum\limits_n
\left(\frac{\mu(S_n)}{|I_n|}\right)^{p-1} \frac{\mu(S_n)}{|I_n|} =
 C \sum\limits_n
\left(\frac{\mu(S_n)}{|I_n|}\right)^{p}. 
$$
We have proved estimate (\ref{83-30}). We conclude that,
by \cite[Theorem 13.1]{gohb},
$$
\|T\|_{{\cal S}_p}^p \le C_1\sum\limits_n
\left(\frac{\mu(S_n)}{|I_n|}\right)^{p}.
$$
Now recall that  $T={\cal J}_1^*{\cal J}_1$ and $p=r/2$. Hence 
$\|{\cal J}_1\|_{{\cal S}_r}^r =\|T\|_{{\cal S}_p}^p<\infty$.
\smallskip
\\
{\bf Proof for the case $0<r\le 1$.} 
Now we prove the statement for $0<r\le 1$ using an idea of \cite{parf1}.
Let $D_n$ be the smallest disc containing the Carleson square $S_n$.
By Lemma 3.5, $dist\, (S_n, \Omega(\Theta, \vep) )\ge |I_n|$,
whence the radius $d_n$ of $D_n$ does not exceed 
$2|I_n|/3$. Let $\tilde D_n$ be a disc with the same center of 
the radius $\tilde d_n = 3|I_n|/4$. It follows that $dist\, (\tilde D_n,
\Omega(\Theta, \vep))\asymp |I_n|$. In this case $\Theta$ 
is analytic in $\tilde D_n$ and $|\Theta(z)|\le C=C(\vep)$, 
$z\in \tilde D_n$. We have $\Theta(z)=1/\overline{\Theta(1/\overline z)}$,
$z\in \tilde D_n\setminus \dd$, and we conclude, analogously to Lemma 5.4,
that 
\begin{equation}
\label{sold}
\|k_z\|_2^2\lesssim |I_n|^{-1}, \qquad z\in \tilde D_n.
\end{equation}

By a well-known Rotfeld's inequality for classes ${\cal S}_r$ 
with $r\le 1$, $\|A+B\|^r_{{\cal S}_r} \le \|A\|^r_{{\cal S}_r}
+\|B\|^r_{{\cal S}_r}$.
Therefore we may represent the embedding operator
${\cal J}_1: {K^2_\Theta} \to L^2(\mu_1)$ as the sum of embedding 
operators $J_n: {K^2_\Theta} \to L^2(\mu_1|_{S_n})$ and estimate their
${\cal S}_r$-norms separately. Now we factorize the operator
$J_n$ as $J_n =J_n^{(2)} J_n^{(1)}$ where $J_n^{(1)}$
is the embedding operator from $K^2_\Theta$ to $H^2(\tilde D_n)$
and $J_n^{(2)}$ is the embedding operator from $H^2(\tilde D_n)$
to $L^2(\mu_1|_{S_n})$. Here $H^2(\tilde D_n)$ denotes
the Hardy space in the disc $\tilde D_n$.

By the standard properties of the classes ${\cal S}_r$, we have
\begin{equation}
\label{eggg}
\|J_n\|_{{\cal S}_r} \le \|J_n^{(2)}\|_{{\cal S}_2}\|J_n^{(1)}\|_{{\cal S}_r}.
\end{equation}
It follows from  (\ref{sold}) that
$$
\|J_n^{(2)}\|^2_{{\cal S}_2} = \frac{1}{2\pi\tilde d_n}
\int_{\partial \tilde D_n}\int_\tz 
|k_z(\zeta)|^2 dm(\zeta)|dz|= 
\frac{1}{2\pi\tilde d_n}
\int_{\partial \tilde D_n} 
\|k_z\|_2^2\, |dz|\le C|I_n|^{-1}.
$$
Note that $d_n \le \delta\tilde d_n$ where $\delta<1$ is an 
absolute constant. Let $s_l$ be the $l$th singular number of 
the operator $J_n^{(1)}$.
Then we have the estimate 
\begin{equation}
\label{sing}
s_l\le C(\delta)\, \delta^{l} \big(\mu(D_n)\big)^{1/2}.
\end{equation}
Indeed, by a translation and linear
change of variables we may assume that $\tilde D_n  =\dd$, 
$D_n=\delta \dd$, and $\nu$ is a measure in $\delta\dd$. Then
$s_l$ does not exceed the norm of the restriction 
of the embedding operator into $L^2(\nu)$ on the subspace $z^l H^2(\dd)$
of $H^2(\dd)$. Note that $|f(z)|\le C(\delta) \|f\|_2$, $z\in \delta\dd$. Then
$$
\|z^l f\|^2_{L^2(\nu)}= \int_{\delta \dd} |z^l f(z)|^2 d\nu(z)\le 
C^2(\delta)\, \delta^{2l} \nu(\delta \dd)\|f\|^2_2, \quad f\in H^2(\dd),
$$
which implies (\ref{sing}). Now, summing $s_l^r$, we obtain 
$\|J_n^{(1)}\|_r^r \le C \big(\mu(S_n)\big)^{r/2}$ 
whence, by (\ref{eggg}),
$$
\|J_n\|_r^r \le C\bigg(\frac{\mu(S_n) }{|I_n|}\bigg)^{r/2}.
$$
We conclude that
$$
\|{\cal J}_1\|_r^r \le C \sum\limits_n 
\bigg( \frac{\mu(S_n) }{|I_n|} \bigg)^{r/2}.
$$

Finally, the case $1<r<2$ follows by interpolation between ${\cal S}_1$
and ${\cal S}_2$ (see \cite[Section 2]{parf1}). 
The proof of Theorem 5.2 is completed. $\bigcirc$

\bigskip

\begin{center}
{\bf 6. Necessary conditions for ${\cal J}_\mu \in {\cal S}_r$. 
Proof of Theorem 1.4}
\end{center}

In this section we consider conditions which are necessary for the inclusion 
${\cal J}_\mu \in {\cal S}_r$, $r\ge 1$.
We will use a general approach from Parfenov's paper \cite{parf1}.
Let $X$ be a Hilbert space of analytic functions in a domain $D$
with the reproducing kernel $K$. Let $\{D_n\}$ be a partition
of $D$ and assume that for any $n$ there exists $w_n\in D_n$
such that, for $z\in D_n$,
\begin{equation}
\label{duda}
|K(z,w_n)|^2 \ge c K(z,z)K(w_n, w_n)
\end{equation}
where $c$ is a positive constant. Consider the discrete measure 
$\nu=\sum_n (K(w_n,w_n))^{-1} \delta_{w_n}$. 
Put $j_n=\Big(\int_{D_n} K(z,z) d\mu(z)\Big)^{1/2}$.
If the embedding operator 
of $X$ into $L^2(\nu)$ is bounded and the embedding
operator of $X$ into $L^2(\mu)$ is in ${\cal S}_r$, $r\ge 1$, then
$\{j_n\} \in \ell^r$ \cite[Theorem 3]{parf1}.  
\medskip
\\
{\bf Proof of Theorem 1.3.} We fix a numeration $R_n$, $n\in\na$,
of the set of the squares $R_{l,m}$ such that
$R_{l,m}\cap \Omega(\Theta, \varepsilon) \ne\emptyset$.
In each of the squares $R_n$ we choose a point $w_n$ 
with $|\Theta(w_n)|<\varepsilon$. If $R_n = R_{l,m}$, 
put $d_n=2^{-l}$. It is easy to show that there exists
a $\delta<1$ which depends only on $\vep$ such that 
$|\Theta(z)|< \delta$, $z \in R_n$. In other words,
$R_n\subset \Omega(\Theta,\delta)$.

We have $\|k_z\|_2^2 \asymp d_n^{-1}$, $z\in R_n$,
with the constants depending only on $\delta$.
Therefore, the sets $D_n = R_n$ and the points $w_n$
satisfy (\ref{duda}). Also, if we put  $\nu= \sum_n d_n^{-1}\delta_{w_n}$,
then it follows from the construction of Carleson curves 
\cite[Chapter VIII, \S 5]{ga} that $\nu$ is a Carleson measure.
Now, let $\mu$ be a measure on $\bigcup\limits_n R_n$,
let ${\cal J}_\mu: K^2_\Theta \to L^2(\mu)$ be the embedding operator. 
Since ${\cal J}_\mu \in {\cal S}_r$, $r\ge 1$, we have $\{j_n\} \in \ell^r$, 
by \cite[Theorem 3]{parf1},  
where 
$j_n=\Big(\int_{R_n} \|k_z\|^2_2 d\mu(z)\Big)^{1/2} \asymp 
\big(\mu(R_n)/d_n\big)^{1/2}$. $\bigcirc$
\medskip

We conclude this section with the proof of Theorem 1.4.
We start with an elementary estimate for inner functions.
\medskip
\\
{\bf Lemma 6.1.}
{\it Let $\zeta \in \mathbb{T}\setminus \sigma(\Theta)$, 
$z\in \mathbb{D}$, and let $|z-\zeta| < A\, 
dist\,(\zeta, \sigma(\Theta))$ for some constant $A\in (0,1)$. 
Then there is a constant $C=C(A)>0$ such that
$$
\log|\Theta(z)| \le - C (1-|z|) |\Theta'(\zeta)|.
$$}
{\bf Proof.}
By the Frostman theorem,
$\Theta_\alpha=\frac{\Theta-\alpha}
{1-\overline\alpha\Theta}$ is a Blaschke product for almost
all $\alpha$ with $|\alpha|<1$, and
$\|\Theta_\alpha-\Theta\|_\infty\to 0$ when $\alpha\to 0$.
We have also $|\Theta'_\alpha(\zeta)|\to |\Theta'(\zeta)|$,
$\alpha\to 0$, whenever
$\zeta\in \tz\setminus \sigma(\Theta)$.
Thus, it suffices to prove the estimate for the case when
$\Theta$ is a Blaschke product. 

Let $B$ be a Blaschke product with zeros $z_n$ and let $z\in\dd$. Then
$$
\log|B(z)|^2=\sum_n\log\left(
1-\frac{(1-|z_n|^2)(1-|z|^2)}{|1-\overline z_n z|^2}
\right).
$$
Also, recall that 
$|B'(\zeta)|=\sum_{n}\frac{1-|z_n|^2}{|\zeta-z_n|^2}$, 
$\zeta\in\tz$. Since $|z-\zeta| < A\,
dist\,(\zeta, \sigma(\Theta))$, we have $|z-\zeta| < A |z_n - \zeta|$
for any $n$. Therefore
$$
(1-A)|\zeta - z_n| < |1-\overline z_n z|< (1+A)|\zeta - z_n|.
$$
Since $\log(1-t) < -t $, $t\in (0,1)$, we have
$$
\log|B(z)|^2 < - \sum\limits_n
\frac{(1-|z_n|^2)(1-|z|^2)}{|1-\overline z_n z|^2}
$$
$$
<- C(A) (1-|z|) \sum\limits_n
\frac{1-|z_n|^2}{|\zeta- z_n|^2} =-C(A)(1-|z|) |B'(\zeta)|. \qquad \bigcirc 
$$

In the next lemma $\{I_n\}$ denotes the system of arcs 
from Lemma 3.5.
\medskip
\\
{\bf Lemma 6.2.}
{\it Let $\Theta\in CLS$. There exists a $\delta \in (0,1)$ such that 
$|\Theta(z)|\le \delta$ for $z=(1-|I_n|/(2\pi)) \zeta$, $\zeta \in I_n$
$($that is, for $z$ on the interior side of the square $S(I_n)$$)$.
Also we have $k_z(z)=\|k_z\|_2^2 \asymp |I_n|^{-1}$, $z\in S(I_n)$. }
\smallskip
\\
{\bf Proof.}
First we show 
that there exist constants $C_j=C_j(\Theta, \varepsilon)>0$, $j=1,2$,
such that
\begin{equation}
\label{83-1}
 C_1|I_n|^{-1} \le |\Theta'(\zeta)| \le C_2|I_n|^{-1}, \qquad \zeta\in I_n.
\end{equation}
We need to prove only the first inequality; the second follows from Lemma 5.4.

By Lemma 3.5, there exists $w\in\Omega(\Theta,\varepsilon)$ such that
$|\zeta-w|\le C_3|I_n|$,
$\zeta\in I_n$, for an absolute constant $C_3>0$. Hence,
$$
|k_w(\zeta)| \ge \frac{1-|\Theta(w)|}{|\zeta- w|}
\ge C_3^{-1} (1-\varepsilon)|I_n|^{-1}.
$$
On the other hand, by an inequality due to Aleksandrov \cite{al2}, 
for a function $\Theta\in CLS$ we have
\begin{equation}
\label{83-10}
|k_w(\zeta)| \le C_4|\Theta'(\zeta)|, 
\qquad w\in\mathbb{D},\  \zeta\in \mathbb{T},
\end{equation}
which implies (\ref{83-1}).
\smallskip

Now fix $\zeta\in I_n$ and put 
$z= (1-|I_n|/(2\pi)) \zeta$.
Since $dist (I_n, \sigma(\Theta)) \ge |I_n|$, we have
$|\zeta-z|<A\, dist (\zeta, \sigma(\Theta))$ for some $A<1$. 
It follows from (\ref{83-1}) and Lemma 6.1 that 
$|\Theta(z)| \le \delta = \exp(-C(A)C_1/(2\pi))$. 
We conclude that $k_z(z)\asymp |I_n|^{-1}$, $z=(1-|I_n|/(2\pi)) \zeta$, 
$\zeta \in I_n$. The estimate $k_z(z)=\|k_z\|_2^2 \asymp |I_n|^{-1}$,
$z\in S(I_n)$ follows from Lemma 3.2. $\bigcirc$
\medskip
\\
{\bf Corollary 6.3.} 
{\it Let $\Theta \in CLS$, $\varepsilon\in (0,1)$. Then there exists 
$\delta \in (0,1)$ such that $|\Theta(z)|\le \delta$, $z\in G$.}
\smallskip
\\
{\bf Proof.} Recall that $G=\mathbb{D}\setminus \bigcup\limits_n S(I_n)$.
We show that there exists $\delta\in (0,1)$ such that 
$|\Theta(z)|\le \delta$, $z\in \partial G\cap \dd$. Since $\partial G$
is a rectifiable Jordan curve, $|\Theta(z)|\le 1$ in $G$ and 
$\partial G\cap\tz$ is of zero Lebesgue measure, it follows that
$|\Theta(z)|\le \delta$, $z\in G$.

Now let $z\in \partial G\cap \dd$. Then there are two possibilities:
either $z=(1-|I_n|/(2\pi)) \zeta$, $\zeta \in I_n$ for some $n$ 
($z$ is on the 
interior side of some square) or there exist two adjacent squares
$S(I_n)$ and $S(I_m)$ with $|I_n|\le |I_m|$ such that
$z=r \zeta$, $1-|I_m|/(2\pi) \le r\le 1-|I_n|/(2\pi)$. Here $\zeta$
is the common endpoint of the arcs $I_n$ and $I_m$.
In the first case $|\Theta(z)|\le \delta_1<1$ by Lemma 6.2.
Note that, by (\ref{83-1}), $|I_n|\asymp |I_m| 
\asymp |\Theta'(\zeta)|^{-1}$. Hence in the second case 
$|\Theta(z)| \le \delta_2<1$ by Lemma 6.1.
$\bigcirc$
\medskip
\\
{\bf Proof of Theorem 1.4.} We start with the sufficiency
of (\ref{23-a}) and (\ref{dya2}).
As before, put $F= \bigcup\limits_n S(I_n)$, 
$G=\mathbb{D}\setminus F$.
It follows from Theorem 1.2 that the embedding operator of
$K_\Theta^2$ into $L^2(\mu|_F)$ is in ${\cal S}_r$.
Now let $R_{n,m}$ be a dyadic square such that 
$R_{n,m} \cap G\ne \emptyset$. Then, by Corollary 6.3, there is a constant
$\delta<1$ such that $R_{n,m} \cap \Omega(\Theta,\delta) \ne \emptyset$.
By (\ref{dya2}), 
$$
\sum\limits_{R_{n,m}\cap G \ne \emptyset} 
(2^n \mu(R_{n,m}))^{r/2}<\infty,
$$
and the inclusion ${\cal J}_{\mu|_G} \in {\cal S}_r$ follows from the 
Luecking theorem.

By Theorem 1.3, condition (\ref{dya2}) is 
necessary for the inclusion ${\cal J}_\mu 
\in {\cal S}_r$ even for general inner functions.
To prove the necessity of (\ref{23-a}), we 
verify the conditions of Parfenov's theorem
for $D_n = S(I_n)$.
By Lemma 6.2, there exist $w_n\in S(I_n)$ 
such that $\|k_{w_n}\|_2^{2} \asymp |I_n|^{-1}$
and 
$$
|k_{w_n}(z)|^2=\left|\frac{1-\Theta(z)\overline{\Theta(w_n)}}
{1-z\overline w_n}\right|^2 \ge C_1|I_n|^{-2}\ge C_2 k_z(z) k_{w_n}(w_n).
$$ 
We used the estimates $|I_n|\lesssim |1-z\overline w_n|$ 
and $k_z(z) \lesssim |I_n|$ (see Lemma 5.4). 
We have $k_{w_n}(w_n)\asymp |I_n|$ and the measure  
$\nu=\sum_n |I_n|\,\delta_{w_n}$ 
is in the class ${\cal C}_2(\Theta)$ by Theorem 3.1, $(i)$.

If $J\in {\cal S}_r$, $r\ge 1$, then, by \cite[Theorem 3]{parf1}, 
we have $\{j_n\}\in \ell^r$, 
$$
j_n=\Big(\int_{D_n} k_z(z) d\mu(z)\Big)^{1/2}.
$$ 
It remains to note that $j_n \asymp \big(\mu(S(I_n))/|I_n| \big)^{1/2}$
since $k_z(z)\asymp |I_n|^{-1}$, $z\in S(I_n)$.
$\bigcirc$
\medskip
\\
{\bf Remarks.} 1. Condition $R_{n,m} \in {\cal R}(\varepsilon, A)$ 
in Theorem 5.2 means that
the distance from the dyadic square 
$R_{n,m}$ to the level set $\Omega(\Theta,\varepsilon)$
is not much larger than its size. The constant $A=12\pi +2$, which appears in
the proof of Theorem 5.2, 
is, by no means, sharp. There is a certain gap between 
the sufficient condition (\ref{m1.8}) and the necessary 
condition (\ref{dya2}). We point out that the inclusion
$R_{n,m} \in {\cal R}(\varepsilon, A)$ does not necessarily imply that
$R_{n,m}\cap\Omega(\Theta,\varepsilon_1) \ne \emptyset$ for some
$\varepsilon_1\in (0,1)$ independent of $n$ and $m$.
\smallskip

2. Theorem 1.4 (or, to be presice, its analog for the upper half-plane)
extends a theorem of Parfenov on embeddings of the Paley--Wiener 
spaces \cite{parf2}:
{\it if $\mu$ is a measure on 
$\rl$ and ${\cal J} :PW_a^2 \to L^2(\mu)$, ${\cal J} f=f$,
then 
${\cal J} \in {\cal S}_r$, $r>0$, if and only if }
$$
\sum\limits_{n\in\mathbb{Z}} \big(\mu([n,n+1])\big)^{r/2}<\infty.
$$
Note that, for the function $\theta(z)=\exp(iaz)$ in $\cp$
(which is a $CLS$ inner function in $\cp$),
the intervals $J_n = [n,n+1]$ have the same properties
as the arcs $I_n$ in Lemma 3.5; namely that, for any 
$\varepsilon\in(0,1)$, $dist \big(J_n, \Omega(\theta, \varepsilon)\big)$
is comparable with the length  of the interval $J_n$.

\bigskip
\noindent
Department of Mathematics and Mechanics, \\
St. Petersburg State University, \\
28, Universitetskii pr., St. Petersburg, \\
198504, RUSSIA
\medskip
\\
E-mail: a.baranov@ev13934.spb.edu 

\end{document}